\documentclass[oneside, a4paper,11pt,reqno]{amsart}
\usepackage{color}

\usepackage{amssymb,amsmath,amsthm}
\usepackage{graphicx}
\usepackage[T1]{fontenc}

\newtheorem{hypo1}{Hypothesis}[section]

\newtheorem{prop1}[hypo1]{Proposition}

\newtheorem{thm1}[hypo1]{Theorem}

\newtheorem{lem1}[hypo1]{Lemma}

\def\I{\mathcal{I}}

\def\F{\mathcal{F}}

\def\PP{\mathbb{P}}

\def\ZZ{\mathbb{Z}}

\def\EE{\mathbb{E}}

\newcommand {\floor}[1] {\left\lfloor {#1} \right\rfloor}

\newcommand{\gap}{g}

\title[A combinatorial approach to a model of constrained random walkers]
       {A combinatorial approach to a model of constrained random walkers}
\author{Thibault Espinasse} 
\address{Institut Camille Jordan, CNRS UMR 5208, Universit\'e de Lyon, Universit\'e Lyon 1, 43, Boulevard du 11 novembre 1918, 69622 Villeurbanne, France.}
\email{Thibault.Espinasse@math.univ-lyon1.fr}
\author{Nadine Guillotin-Plantard} 
\email{nadine.guillotin@univ-lyon1.fr}
\author{Philippe Nadeau}
\email{philippe.nadeau@math.univ-lyon1.fr}


\begin{document}

\begin{abstract} 
In \cite{BCEN}, the authors consider a random walk $(Z_{n,1},\ldots,Z_{n,K+1})\in \ZZ^{K+1}$ with the
constraint that each coordinate of the walk is at distance one from the
following one. A functional central limit theorem for the first coordinate is proved and the limit variance is explicited.
In this paper, we study an extended version of this model by conditioning the extremal coordinates to be at some fixed distance at every time. We prove a functional central limit theorem for this random walk. Using combinatorial tools, we give a precise formula of the variance and compare it with the one obtained in \cite{BCEN}.

\end{abstract}

\maketitle

\section{Introduction and results}
Central limit theorems for additive functionals of Markov chains have attracted continuing interest for over half a century. One approach, due to Gordin  \cite{Gor},  rests on martingale approximation. Roughly speaking in good cases the asymptotic normality for functionals of Markov chain is derived from the central limit theorem for martingales. 
Various conditions for the central limit theorem to hold are now known and several expressions for the asymptotic variance can be found in the literature (see for instance \cite{HR}), however it is usually quite difficult to compute it theoretically.
In this note, we are faced with this problem in the study of a generalized version of the so-called "prisoners model" introduced in \cite{BCEN}. We are able to give an explicit formula for the asymptotic variance of the Markov chain we are interested in.
However our approach is not robust enough to give a full description of the asymptotic variance for more general models (see Remark (iii) in Section \ref{sec:conclusion}).

In what follows we will use the notation $[\![1;n]\!] :=\{1,\ldots,n\}$.
Let $K$ be a positive integer. For any $h\in[\![0;K]\!]$, we define
 $$ \mathcal{C}_{K, h}=\{ z\in \ZZ^{K+1}; \forall i\in [\![1;K]\!], |z_{i+1}-z_i| =1  \mbox{  and  }z_{K+1}- z_1= h\}.$$ 
The set $ \mathcal{C}_{K, h}$ is empty unless $K-h$ is even; so let us set $K-h=2\gap$. The set $\mathcal{C}_{K,0}$ corresponds to the Bernoulli bridges with length $K$. 

We can define a neighbourhood structure on $\mathcal{C}_{K,h}$ through the following relation
$$\forall z,z' \in \mathcal{C}_{K,h}, z \sim z' \Leftrightarrow \forall i \in [\![1;K+1]\!], |z_i-z'_i| = 1.  $$
The set of neighbors of $z\in\mathcal{C}_{K,h}$ will be denoted by $\Gamma(z)$ and its cardinality by $\deg_{\mathcal{C}_{K,h}}(z)$.

Let us denote by $(Z_n)_{n\geq 0}$  the Markov chain defined on  $\mathcal{C}_{K,h}$ corresponding to $K+1$
 simple random walks on $\ZZ$ under the shape constraint. In other words, $(Z_n)_{n\geq 0}$ is the Markov chain with state space 
 $\mathcal{C}_{K,h}$ and transition probabilities given by 
$$\PP[Z_{n+1} = z' | Z_{n} = z] = \left\{ 
\begin{array}{ll}
\frac{1}{\deg_{\mathcal{C}_{K,h}}(z)} & \text{if} \  z'\sim z,\\
\  \   0    & \text{otherwise.}
 \end{array}
 \right.
 $$
 
We will denote by $Z_{n,i}$ the $i^{\text{th}}$ coordinate of $Z_n$ and by $\lfloor x \rfloor$ the integer part of a real number $x$. We assume that $Z_{0,1}=0$ almost surely.

In this paper we are interested in the distributional limit of the first coordinate $(Z_{n,1})$  as $n$ tends to infinity. Actually, the convergence to the Brownian motion is not surprising, and we are mainly interested here in the exact value of the variance. 

\begin{thm1}\label{Th1}
The sequence of random processes $\left(\frac{Z_{\lfloor nt\rfloor,1}}{\sqrt{n}}\right)_{t\geq 0, n\geq 1}$  weakly converges to the Brownian motion with variance 
$$\sigma_{K,h}^2=\frac{1}{K} \frac{A_{K,h}}{B_{K,h}}$$
where 
$$A_{K,h}= \sum_{k=0}^{\floor{K/2}}  \binom{K}{2k+1} \binom{K-2k}{\gap-k}$$
and $$B_{K,h}=\sum_{k=0}^{\floor{K/2}} \binom{K}{2k}\  \binom {K-2k}{\gap-k}.$$
The variance $\sigma_{K,h}^2$ has the following properties:
\begin{itemize}
\item[i)] For any $h\in [\![0;K]\!]$, as $K$ tends to infinity,
$$\sigma_{K,h}^2 \sim \frac{2}{K}.$$
\item[ii)] For any $K$,  
$$\sigma_{K,0}^2 <  2 / K .$$ 
\item[iii)] For $K$ large enough, 
$$\sigma_{K,0}^2 > 2 / (K+2).$$
\end{itemize}
\end{thm1}
With our notations, the variance obtained in \cite{BCEN} is equal to $2/ (K+2)$, so item iii) means that for large $K$,
our variance is larger than the one obtained in the unconstrained case. We conjecture that this should hold for any $K$.

The paper is organized as follows: in Section~\ref{sec:comb}, various enumerative results about walks are proved. They will be used in Section~\ref{sec:proba}. 
Theorem \ref{Th1} is proved in Section~\ref{sec:proba}. Our approach is based on a decomposition of $Z_{n,1}$ as a sum of a martingale and a bounded function. 
We will compute explicitly the two parts, and provide a geometrical interpretation for the decomposition. The asymptotic properties of the variance will be proved using local limit theorems for random walks on $\ZZ$. 
In Section~\ref{sec:unconstrained} we explain how to adapt our method to derive the main result of \cite{BCEN}.  
In Section~\ref{sec:conclusion} some variations of the model are proposed and open problems are discussed. 

\section{Combinatorics}
\label{sec:comb}

We prove various enumerative results that will be used in Section
 \ref{sec:proba}.

If $z\in\mathcal{C}_{K,h}$, its shape is defined as $\F(z)=(F_1,\ldots,F_K)$ by $F_i=z_{i+1}-z_{i}$. We refer naturally to the elements of $\F(z)$ as the {\em steps} of $z$, which can be either up or down. The set of possible shapes is given by 
\[Sh_{K,h}=\{(F_i)_i\in\{\pm 1\}^{K}; \sum_{i=1}^KF_i=h\}.\]
 Clearly $z\in\mathcal{C}_{K,h}$ if and only if $(z_1,\F(z))\in\ZZ\times Sh_{K,h}$. Let us denote, for every $z\in\mathcal{C}_{K,h}$, 
$$ \Gamma^{\pm }(z)=\{ z'\sim z; z_1'-z_1= z_{K+1}'-z_{K+1} = \pm 1\}.$$

\begin{figure}[!ht]
\includegraphics{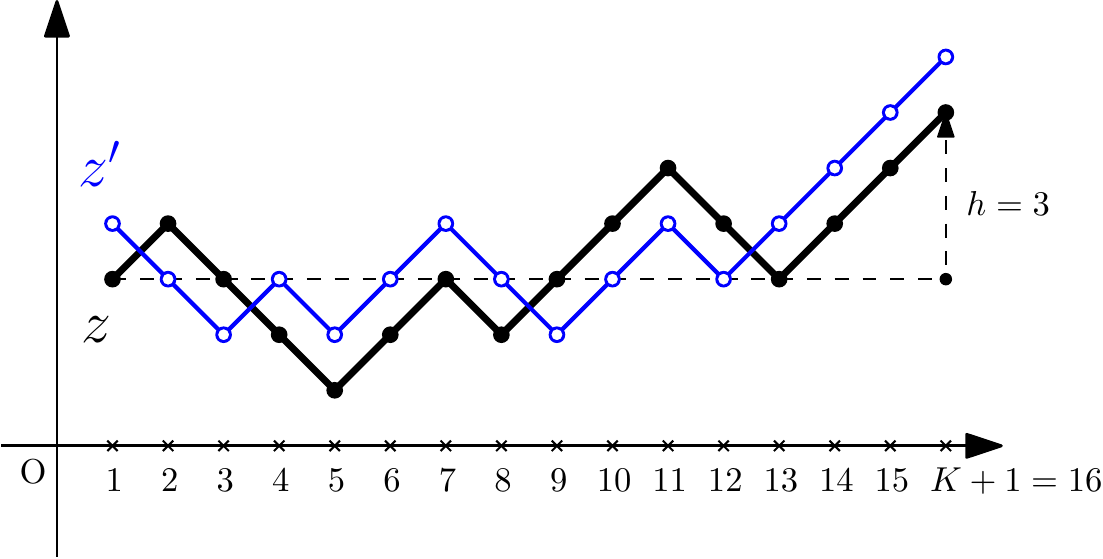}
\caption{Two paths $z\sim z'$ with $z'\in\Gamma^+(z)$
\label{Figure_Chemins}}
\end{figure}

Now let $z,z'\in\mathcal{C}_{K,h}$ be such that $z'\in\Gamma^+(z)$, and let $F=\F(z),F'=\F(z')$ be the shapes of $z,z'$ respectively. For any $i\in \{1,\ldots,K\}$, define $M_i=\pm 1$ if $F_i=F'_i=\pm 1$ and $M_i=0$ otherwise. So $M=(M_1,\ldots,M_{K})$ has values in $\{1,0,-1\}$, and can thus be represented as a Grand Motzkin path; we define $\Phi^+(z,z')=M$.
 Note that $\sum_i M_i=h$, so this path $M$ goes from $(0,0)$ to $(K,h)$; we note $\mathcal{M}_{K,h}$ the set of such sequences, which has cardinality
\begin{equation}\label{GMP}
|\mathcal{M}_{K,h}|=\sum_{k=0}^{\floor{K/2}}\binom{K}{2k}\binom{K-2k}{\gap-k}.
\end{equation}
 Now if we fix the starting point of $z$, then from any element of $\mathcal{M}_{K,h}$ we can clearly reconstruct a pair $z,z'$, so we have the following result

\begin{prop1}
\label{prop:bijection_motzkin}
For any $z_1$, $\Phi^+$ is a bijection between $\{(z,z')\in\mathcal{C}_{K,h}^2; z'\in\Gamma^+(z)\}$ and $\mathcal{M}_{K,h}$.
\end{prop1}

Denote by $A(z)$ the algebraic area between the path $z$ and the $X$-axis. Note that 
\[A(z) = \frac{1}{2}z_{1}+ z_2 + \cdots + z_{K} + \frac{1}{2}z_{K+1}\]
 
\begin{figure}[!ht]
\includegraphics{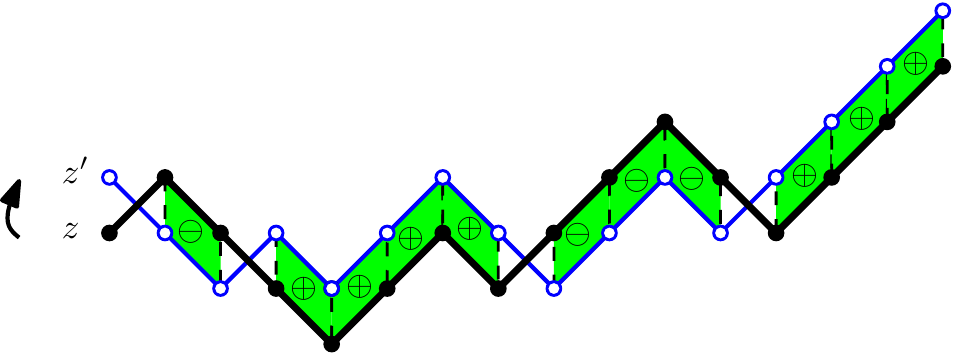}
\caption{\label{Figure_Area}
The  area $A(z')-A(z)$ is equal to $-1+4-3+3=3$}
\end{figure}

We now prove that for any fixed path $z\in \mathcal{C}_{K,h} $ the sum over $z'\in\Gamma(z)$ of the algebraic areas boils down to zero. Fix a path $z\in \mathcal{C}_{K,h} $. Let $(z',p)$ be a {\em marked path} (with respect to $z$), which we define as a path $z'\in\Gamma(z)$ together with a step $p\in z'$ which does not cross $z$.  Define 
 $sign(z;z',p)=1$ (\emph{resp. $-1$}) if the step $p$ occurs above (\emph{resp.} below) $z$. With these notations it is clear that $A(z')-A(z)=\sum_p sign(z;z',p)$, where the sum is over all marking steps in $z'$; see Figure~\ref{Figure_Area}.

\begin{figure}[!ht]
\includegraphics[width=0.9\textwidth]{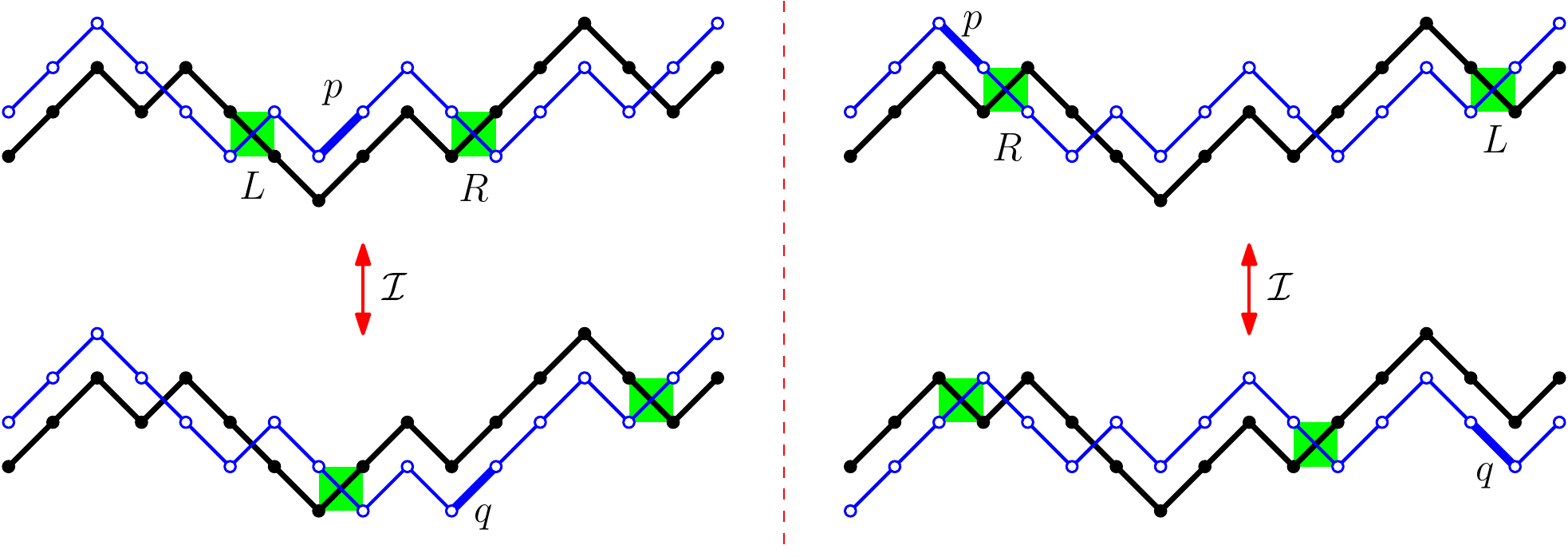}
\caption{Two cases of the application of the involution $\I$.
\label{Figure_Involution}}
\end{figure}

We now define an involution $\I$ of these marked paths. Note first that $z'$ crosses $z$ an even number of times, and let $2k$ be this quantity. The case $k=0$ occurs only when $z'\in\{z^+,z^-\}$ which are defined as the two paths obtained from $z$ by shifting it one unit up or down. In this case define simply $\I(z;z^\pm,p)=(z;z^\mp,p)$. Suppose now $k\neq 0$, and denote by $L$ (\emph{resp.} $R$) the nearest crossing left (\emph{resp.} right) of $p$. There are two special cases: if $p$ occurs before the first crossing of $z'$, then $L$ is defined as the last crossing of $z'$, while if $p$ occurs after the last crossing of $z'$, then $R$ is defined as the first crossing of $z'$. Examples of the definition of $L,R$ are illustrated in Figure~\ref{Figure_Involution}.

Now exactly one step among $\{L,R\}$ is a step of $z'$ of a different up/down type as $p$. Let $q$ be this crossing step, and exchange $p$ with $q$ in $z'$, keeping all other steps unchanged. This defines a new path $z''\in \Gamma(z)$, in which $p$ becomes crossing while $q$ is noncrossing. So $(z'',q)$ is a marked path, and we define  $\I(z;z',p):=(z;z'',q)$. We have the following lemma, whose proof is immediate.

\begin{lem1}
For any $z$, $\I$ is an involution on marked paths which is sign reversing.
\end{lem1}
  
\begin{prop1}
\label{prop:zerosum}
For any $z\in\mathcal{C}_{K,h}$, $\displaystyle{ \sum_{z'\in \Gamma(z)}  \left(A(z') - A(z)\right) = 0.}$
\end{prop1}

\begin{proof}
The sum to compute is equal to the sum of $sign(z;z',p)$ for all marked paths $(z',p)$. By the previous lemma, the involution $\I$ pairs such marked paths two by two with opposite signs, so the sum boils down to zero.
\end{proof}

\begin{prop1}
\label{prop:total_area}
For any $K,h$, $\displaystyle{\sum_{\stackrel{z,z'\in\mathcal{C}_{K,h}}{z'\in\Gamma^+(z)}}[ A(z') - A(z)]=\sum_{k=0}^{\floor{K/2}}  \binom{K}{2k+1} \binom{K-2k}{\gap-k}}.$
\end{prop1}

\begin{proof}
The l.h.s. can be rewritten as the sum of $sign(z;z',p)$ over all $z,z',p$ such that $z'\in\Gamma^+(z)$ and $(z',p)$ is marked. Define $(z,z'',q):=\I(z,z',p)$. If $z''$ also belongs to $\Gamma^+(z)$, then the corresponding terms in the sum will cancel each other. The surviving terms correspond to $z''\in\Gamma^-(z)$, which by inspection occurs precisely in the following cases:\begin{enumerate}
\item \label{it1} $p$ occurs before the first crossing, and $p$ is a descent; 
\item \label{it2} $p$ occurs after the last crossing, and $p$ is an ascent; 
\item \label{it3} $z'=z^+$ and $p$ is any step.
\end{enumerate}

In all three cases, $sign(z,z',p)$ is always equal to $1$. Thus we need to enumerate the triplets $(z;z',p)$ verifying \eqref{it1}, \eqref{it2} or \eqref{it3}. 

There are clearly $K \binom{K}{\gap}$ triplets satisfying \eqref{it3}. We now compute the number of triplets $(z,z',p)$ verifying case \eqref{it1} or  \eqref{it2}, and such that there are $2k>0$ crossings between $z$ and $z'$. The position of $p$ together with the position of these crossings form a subset $\{j_1<j_2<\cdots<j_{2k+1}\}$ of $\{1,\ldots,K\}$ ; here the position of $p$ is $j_1$ in case \eqref{it1} and $j_{2k+1}$ in case \eqref{it2}. There are $\binom{K}{2k+1}$ such subsets, and, conversely, the knowledge of such a subset determines the positions of $p$ and the crossings of $z$ and $z'$ in each case. The remaining $K-2k-1$ steps are the same in $z$ and $z'$, so we only need to count possible steps for $z$. We want to ensure that $z$ belongs to $\mathcal{C}_{K,h}$. Notice that among the crossings there are $k$ up steps and $k$ down steps. Now in case  \eqref{it1}, $p$ is a down step, so there must be exactly $g-k-1$ down steps among the remaining $K-2k-1$ steps; in case \eqref{it2}, $p$ is an up 
step, so there must 
be exactly $g-k$ down steps among the remaining $K-2k-1$ steps. In total, this represents $\binom{K-2k-1}{g-k-1}+\binom{K-2k-1}{g-k}=\binom{K-2k}{g-k}$ possibilities.

Adding everything up we obtain the desired expression (note that the case $k=0$ corresponds to \eqref{it3}). \end{proof}

%

\section{Proof of Theorem \ref{Th1}}
\label{sec:proba}

\subsection{Martingale}
Denote by $(\mathcal{F}_n)_{n\geq 0}$  the natural filtration of the Markov chain $Z_n$.

\begin{lem1}\label{Le1}
The sequence $(A (Z_n))_{n\geq 0}$ 
is a martingale with respect to the filtration $(\mathcal{F}_n)_{n\geq 0}$.
\end{lem1}
\begin{proof}
Since $(Z_n)_n$ is a Markov chain, it is enough to prove that for any $z\in\mathcal{C}_{K,h}$, 
$$\mathbb{E}[ A(Z_{n+1}) - A(Z_n)   | Z_n=z]=0$$
or equivalently (since the Markov chain moves uniformly on its neighbours) for any $z\in\mathcal{C}_{K,h}$, 
$$ \sum_{z'\sim z}  [ A(z') - A(z)] = 0. $$
We recognize here the content of Proposition~\ref{prop:zerosum}.
\end{proof}

\subsection{The coupling}
\begin{figure}[!ht]
\includegraphics[width=0.2\textwidth]{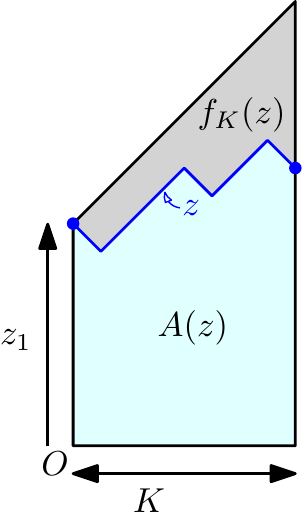}
\caption{Relation between the first coordinate $z_1$ and the area under the path.
\label{Figure_Coupling}}
\end{figure}

Denote by $f_K(z)$ the area between a path $z=(z_1,z_2,\ldots, z_{K},z_{K+1})\in  \mathcal{C}_{K,h}$,
 and the segment $(z_1, z_1+1, ...., z_1+K)$.
Then, for any $z\in  \mathcal{C}_{K,h}$ (see Figure~\ref{Figure_Coupling}),
\begin{equation}\label{aire}
f_K(z) +A(z) = K z_1 + \frac{K^2}{2}.
\end{equation}

Since $f_{K}$ is bounded by $K^2$, the position of the first walker and the area under the path can be asymptotically related as 
\begin{equation}\label{aire2}
\frac{Z_{[nt],1}}{\sqrt{n}} = \frac{1}{K}\frac{ A (Z_{[nt]})}{\sqrt{n}}+o(1).
\end{equation}
Theorem \ref{Th1} will then be deduced from the functional central limit theorem for the rescaled martingale $(A(Z_{n})/\sqrt{n})_{n\geq 1}.$

If $z\in\mathcal{C}_{K,h}$, its shape is defined as $\F(z)=(F_1,\ldots,F_{K+1})$ by $F_i=z_{i}-z_{1}$. The set of possible shapes is given by
 $$\overline{Sh}_{K,h}=\{z\in\mathcal{C}_{K,h}; z_1=0\}.$$ 
This description of the possible shapes of $\mathcal{C}_{K,h}$ differs from the one used in Section \ref{sec:comb} but will be more convenient in the next computations.
Note that a natural neighbourhood structure is given on $\overline{Sh}_{K,h}$ by the relation
$$F \sim_S F' \Leftrightarrow \forall i \in [\![1;K+1]\!], F_i - F'_i \in \{0,2\} \text{ or } \forall i \in [\![1;K+1]\!], F_i - F'_i \in \{0,-2\} . $$
If $F,F'$ are neighbours in $\overline{Sh}_{K,h}$, then one of following statements is true:
\begin{itemize}
 \item If $F_i - F'_i \in \{0,2\}$ then $F-(1, \cdots, 1)$ and $F'$ are neighbours in $\mathcal{C}_{K,h}$. 
\item If $F_i - F'_i \in \{0,-2\}$ then $F+(1, \cdots, 1)$ and $F'$ are neighbours in $\mathcal{C}_{K,h}$.
\end{itemize}
Note that the special case $F = F'$ gives two neighbours in $\mathcal{C}_{K,h}$ while the others cases provide a one-to-one transformation.\\
Conversely, if $z,z'$ are neighbours in $\mathcal{C}_{K,h}$, then either $\F(z) = \F(z')$ and $z = z' \pm (1, \cdots, 1)$, or $\F(z)$ and $\F(z')$ are distinct neighbours in $\overline{Sh}_{K,h}$. 

\noindent We denote by $\deg_S(F)$ the number of neighbours of $F\in \overline{Sh}_{K,h}$, i.e.
$$\deg_S(F) := \sharp \left\{G \in \overline{Sh}_{K,h}; G \sim_S F \right\}. $$
According to the previous remark, we have
$$\deg_{\mathcal{C}_{K,h}}(z) = \deg_S(\F(z))+1. $$
Consider the Markov chain $(F_n)_{n}$ with values in $\overline{Sh}_{K,h}$, with transition probabilities given by
$$\displaystyle p(F,G) = \left\{ \begin{array}{lll}
\frac{1}{ \deg_S(F) +1} &  \mbox{  if  } & G \sim_S F \mbox{  and  } G\neq F,
\\ \frac{2}{\deg_S(F)+1} & \mbox{ if } & G=F.
\end{array}
\right.$$
This irreducible and ergodic Markov chain is reversible with stationary measure 
$$\pi_{S}(F) := \frac{\deg_S(F) +1} { \sum_{F\in\overline{Sh}_{K,h} } (\deg_S(F) +1)}.$$
The Markov chains $(\F(Z_n))_{n \in \mathbb{N}}$ and $(F_n)_{n \in \mathbb{N}}$ are then identically distributed.
Indeed, 
\begin{align*}\mathbb{P}[\F(Z_{n+1}) = G|Z_n = z]& = \frac{1}{\deg_{\mathcal{C}_{K,h}}(z)}\sum_{z'\sim z} 1_{\F(z') = G} 
\\& =\left\{ 
\begin{array}{ll}
\frac{1}{\deg_{S}(\F(z))+1} & \text{if}  ~ G \sim_S \F(z), G \neq \F(z),
\\ \frac{2}{\deg_{S}(\F(z))+1}    & \text{if} ~ G =\F(z)
 \end{array}
 \right.      
\\ & = \mathbb{P}[F_{n+1} = G|F_n = \F(z)] \end{align*}
Hence, 
\begin{align*}
\mathbb{P}[\F(Z_{n+1}) = G|\F(Z_n) = F]& = \sum_{z \in \mathcal{C}_{K,h}} \mathbb{P}[\F(Z_{n+1}) = G|\F(Z_n) = F; Z_n = z]\mathbb{P}[Z_n = z |\F(Z_n) = F]
\\ & =  \mathbb{P}[F_{n+1} = G|F_n = F]  \Big(\sum_{z \in \mathcal{C}_{K,h}} \mathbb{P}[Z_n = z |\F(Z_n) = F] \Big)
\\ & = \mathbb{P}[F_{n+1} = G|F_n = F].
\end{align*}
Let us remark that for every $z,z' \in \mathcal{C}_{K,h}$ s.t.  $z'\sim z$ with $z'\neq z$,  the difference of area $|A(z')-A(z)|$ can be computed with respect to $\F(z')$ and $\F(z)$ only. 
Indeed, we have the following relations 
\begin{eqnarray*}
|A(z') - A(z)|&=&  \left\{ 
\begin{array}{lll}
 |A(\F(z')) - A(\F(z)+(1,\cdots,1)) | & \text{ if }  \exists i \text{ s.t. } \F(z')_{i} -\F(z)_{i} = 2,\\
|A( \F(z')) - A(\F(z)-(1,\cdots,1) ) |  & \text{ if } \exists i \text{ s.t. } \F(z')_{i} -\F(z)_{i} = -2,\\
  \hspace{2.5cm}  K  &  \text{ if } \F(z')=\F(z).
\end{array}
\right.
\end{eqnarray*}
The quantity  $|A(z')-A(z)|$ will then be denoted by $\delta A(\F(z),\F(z'))$.
From the central limit theorem for martingales (see for instance \cite{Rick}, Theorem 7.4, p.374), we have to compute the a.s. limit as $n$ tends to infinity of
\begin{eqnarray*}
\Sigma_n &:=& \frac{1}{n} \sum_{k=1}^{n} \mathbb{E}[ (A(Z_k) -A(Z_{k-1}) )^2 |Z_{k-1}]\\
&=& \frac{1}{n} \sum_{k=1}^{n} \mathbb{E}[ \delta A(\F(Z_{k-1}),\F(Z_k))^2 |\F(Z_{k-1}); Z_{k-1,1} ]\\
&= & \sum_{F\in \overline{Sh}_{K,h}} \frac{1}{n} \sum_{k=1}^{n} 1_{ \{ \F(Z_{k-1}) = F\}}  \mathbb{E}[ \delta A(\F(Z_{k-1}),\F(Z_k))^2|\F(Z_{k-1}) =F ]\\
&=& \sum_{F\in \overline{Sh}_{K,h}} \Big( \frac{1}{n} \sum_{k=1}^{n} 1_{ \{\F(Z_{k-1}) = F\}} \Big) \mathbb{E}[ \delta A(\F(Z_{0}),\F(Z_1))^2  |\F(Z_0) =F ]\\
&\stackrel{law}{=}& \sum_{F\in \overline{Sh}_{K,h}} \left( \frac{1}{n} \sum_{k=1}^{n} 1_{ \{F_{k-1} = F\}} \right) \mathbb{E}[ \delta A (F_0,F_1 )^2 |F_0 =F ]
\end{eqnarray*}
Since the Markov chain $(F_n)_n$ is ergodic with invariant measure $\pi_{S}$, the sequence of random variables $(\Sigma_n)_n$ converges almost surely to the constant 
$$ \sum_{F\in \overline{Sh}_{K,h}} \mathbb{E} [ (\delta A (F_0,F_1))^2| F_0=F]  \pi_S(F),$$
which can be rewritten as $\mathbb{E}_{\pi_S} [ (A( Z_1) -A (Z_0))^2]$ by remarking that $\F(Z_0)=Z_0$.
Then, using (\ref{aire2}), we obtain that
$$\sigma_{K,h}^2 = \frac{1}{K^2}   \mathbb{E}_{\pi_S} [ (A( Z_1) -A (Z_0) )^2].$$
From (\ref{aire}), we deduce that
$$f_K(Z_1) - f_K(Z_0) = (K-1) (Z_{1,1}- Z_{0,1} )  - [A(Z_1) - A(Z_0)].$$
From the reversibility of the Markov chain $(F_n)_n$, we have 
$$ \mathbb{E}_{\pi_S} [ (A(Z_1) -A(Z_0) ) (f_K(Z_1) -f_K(Z_0) ) ] =0$$
Indeed, $\mathbb{E}_{\pi_S} [ (A(Z_1) -A(Z_0) ) f_K(Z_1) ]$ can be rewritten as
\begin{eqnarray*}
& &\sum_{F\in \overline{Sh}_{K,h}}\pi_S(F) \sum_{F'\neq F} p(F,F') (A(z_1)-A(z_0)) f_K(z_1) {\bf 1}_{\{z_0\sim z_1;\F(z_0)=F;\F(z_1)=F'\}} \\
&=& \sum_{(F,F'); F\neq F'} \pi_S(F') p(F',F)  (A(z_1)-A(z_0))f_K(z_1)  {\bf 1}_{\{z_0\sim z_1;\F(z_0)=F;\F(z_1)=F'\}}  \\
&=& - \sum_{F'} \pi_S (F') \sum_{F\neq F'} p(F',F)  (A(z_0)-A(z_1))f_K(z_1) {\bf 1}_{\{z_0\sim z_1;\F(z_0)=F;\F(z_1)=F'\}} \\
&=& -\mathbb{E}_{\pi_S} [ (A(Z_1) -A(Z_0) ) f_K(Z_0) ] 
 \end{eqnarray*}
and the last expectation is zero since for any fixed shape $F\in \overline{Sh}_{K,h}$, we have 
 \begin{eqnarray*}
 \mathbb{E}[(A(Z_1)-A(Z_0)) f_K(Z_0)|\F(Z_0) =F] &=& \mathbb{E}[(A(Z_1)-A(Z_0)) f_K(Z_0)| Z_0 =F] \\
 &= & f_K(F)  \mathbb{E}[ A(Z_1)-A(z_0) | Z_0 = F] \\
 &=& 0
\end{eqnarray*}
from Lemma \ref{Le1}.
Therefore, 
$$\sigma_{K,h}^2 = \frac{1}{K}\  \mathbb{E}_{\pi_S} [ (A(Z_1) -A(Z_0)) ( Z_{1,1}-Z_{0,1}) ].$$
From Proposition~\ref{prop:zerosum} and the expression of the measure $\pi_{S}$, we deduce 
\begin{eqnarray*}
\sigma_{K,h}^2 &=& \frac{1}{K}\  \frac{A_{K,h}}{B_{K,h}}
\end{eqnarray*}
where
$$A_{K,h}= 2 \sum_{\stackrel{\small z,z'\in\mathcal{C}_{K,h}}{ z'\in\Gamma^+(z)}} (A(z')-A(z))\ 
\text{ and }\ 
B_{K,h}=\sum_{F\in\overline{Sh}_{K,h}} (\deg_S (F)+1).$$
Applying Propositions~\ref{prop:total_area} and \ref{prop:bijection_motzkin} combined with (\ref{GMP}), we get
$$A_{K,h}= 2\sum_{k=0}^{\floor{K/2}}  \binom{K}{2k+1} \binom{K-2k}{\gap-k}\  \text{ and }\ 
B_{K,h}=2|\mathcal{M}_{K,h}|=2 \sum_{k=0}^{\floor{K/2}}\binom{K}{2k}\binom{K-2k}{\gap-k}.$$

\subsection{Properties of the variance}
\noindent We are interested in the properties of the variance $\sigma_{K,h}^2.$
Denote by $(S_n)_{n\geq 0}$ the random walk on $\mathbb{Z}$ starting from 0 and moving according to the following rule
$$\PP[S_{n+1} = y | S_n = x] =\frac{1}{3} \mbox{  if   } y \in\{x-1,x,x+1\}.$$
Then, for any positive integer $K$ and any $h\in [\![0;K]\!]$ with $K-h$ even,
$$\PP[S_{K} = h]=\frac{ |\mathcal{M}_{K,h}| }{3^{K}}$$
where $|\mathcal{M}_{K,h}| $ denotes the number of Grand Motzkin paths from $(0,0)$ to $(K,h)$.
The variance $\sigma_{K,h}^2$ can then be rewritten 
\begin{equation}\label{tll}
\sigma_{K,h}^2 = \frac{1}{K} \frac{\PP[S_{K}=h+1]+\PP[S_{K}=h-1]}{\PP[S_{K}=h]}.
\end{equation}		
Item i) directly follows from the local central limit theorem for the random walk $(S_n)_{n\geq0}$ ( see for instance Proposition P4 p.46 in \cite{Spi}).

Remark that from the symmetry of the random walk $(S_n)_n$, we have
$$\sigma_{K,0}^2 = \frac{2}{K} \frac{\PP[S_{K}=1]}{\PP[S_{K}=0]}.$$
From Fourier inversion formula, for any $x\in \ZZ$, we have
\begin{eqnarray*}
\PP[S_n= x]  &=&\frac{1}{2\pi} \int_{-\pi}^{\pi} e^{i\theta x} \EE\Big[ e^{i\theta S_n}\Big]\, d\theta\\
&=& \frac{1}{2\pi} \int_{-\pi}^{\pi} e^{i\theta x} \EE\left[ e^{i\theta S_1}\right]^n\, d\theta\\
&=&  \frac{1}{2\pi} \int_{-\pi}^{\pi} e^{i\theta x} \left( \frac{1+2\cos  \theta}{3} \right)^n\, d\theta\end{eqnarray*}
since the increments of the random walk $(S_n)_n$ are independent and identically distributed according to the uniform distribution  on the discrete set $\{-1,0,+1\}$.
Therefore, we deduce  
\begin{eqnarray*}
3^{K} \PP[S_{K}=1] &=& (2 \pi )^{-1} \left|\int_{-\pi}^{\pi} ( 1 + 2 \cos(\theta))^{K} e^{i\theta} d\theta\right|\\
&< & (2 \pi )^{-1} \int_{-\pi}^{\pi} ( 1 + 2 \cos(\theta))^{K} d\theta = 3^{K} \PP(S_{K}=0),
\end{eqnarray*}
which yields assertion ii).

Let us now prove item iii).  
Since the increments of the random walk $(S_n)_n$ are symmetric and bounded, from Theorem 2.3.5 (Formula (2.23)) in \cite{LL}, there exists a constant $c$ such that for any  $n$ and $x$, 
\begin{equation}\label{vlada}
\Big| \PP[S_{n} = x] - \frac{\sqrt{3} e^{-3 x^2/(4n)}}{2\sqrt{\pi n}} \Big| \leq \frac{c}{n^{3/2}}.
\end{equation}
Now concerning iii) it is enough to prove that the limit as $m$ tends to infinity of the sequence 
$$ u_K:=\frac{2\ \PP[S_{K}=0]}{(K+2) (\PP[S_{K}=0] -\PP[S_{K}=1])}$$ is strictly greater than one. 
From (\ref{vlada}), we have 
\begin{eqnarray*}
u_K&=&\frac{2}{K+2}\  \frac{ \frac{\sqrt{3}}{\sqrt{4\pi K}}+\mathcal{O}(K^{-3/2})}{  \frac{\sqrt{3}}{\sqrt{4\pi K}}(1- e^{-3/4K})+\mathcal{O}(K^{-3/2})  }\\
&=& \frac{2}{K+2}\  \frac{ \frac{\sqrt{3}}{\sqrt{4\pi K}}+\mathcal{O}(K^{-3/2})}{  \frac{3\sqrt{3}}{4\sqrt{4\pi} K^{3/2}}+\mathcal{O}(K^{-3/2})  }\\
&\sim & \frac{8}{3} \ \text{as}\  K\rightarrow +\infty.
\end{eqnarray*}

\section{The unconstrained case}
\label{sec:unconstrained}

In this section we indicate how to adapt the ideas we developed in order to give an alternative proof of the main result of~\cite{BCEN}. We now consider the Markov chain $(Z_n^{*})_{n\geq 0}$ with values in the state space $\mathcal{C}_{K}=\{ z\in \ZZ^{K+1}; \forall i\in [\![1;K]\!], |z_{i+1}-z_i| =1\}$, so that $\mathcal{C}_{K}$ is the union $\cup_h\mathcal{C}_{K, h}$. The difference is that both ends of the walks are allowed to move in different directions. The set of paths $\Gamma(z)=\Gamma^+(z)\sqcup\Gamma^-(z)$ is defined as before.

 By Proposition~\ref{prop:bijection_motzkin}, $\Phi^+$ is a bijection between $\{(z,z')\in\mathcal{C}_{K}^2; z'\in\Gamma^+(z)\}$ and  $\mathcal{M}_{K}:=\cup_h\mathcal{M}_{K,h}$, which clearly has cardinality $3^K$ since each step can be chosen independently. Now modify the area of $z\in \mathcal{C}_{K}$ by setting
\[A^{*}(z) := \frac{3}{2}z_{1}+ z_2 + \cdots + z_{K} + \frac{3}{2}z_{K+1}=z_1+A(z)+z_{K+1}\]

Extend any path $z\in\mathcal{C}_{K}$ by unit horizontal steps at both ends: then $A^{*}(z)$ is the algebraic area below this extended path. We have the following proposition which is the counterpart of Propositions~\ref{prop:zerosum} and~\ref{prop:total_area} in this unconstrained setting.
\begin{prop1}\label{prop1}
For any $z\in\mathcal{C}_{K}$, $\displaystyle{ \sum_{z'\in \Gamma(z)}  \left(A^{*}(z') - A^{*}(z)\right) = 0.}$

For any $K$, $\displaystyle{\sum_{\stackrel{z,z'\in\mathcal{C}_{K}}{z'\in\Gamma^+(z)}}[ A^{*}(z') - A^{*}(z)]=2\cdot 3^K}$.
\end{prop1}

\begin{proof}[Sketch of the proof]
The difference $A^{*}(z')-A^{*}(z)$ between two paths can be written as the sum of $sign(z;z',p)$ over all marked paths $(z;z',p)$ where the mark $p$ can now also be one of the two extra horizontal steps, see figure below.
\begin{center}
\includegraphics{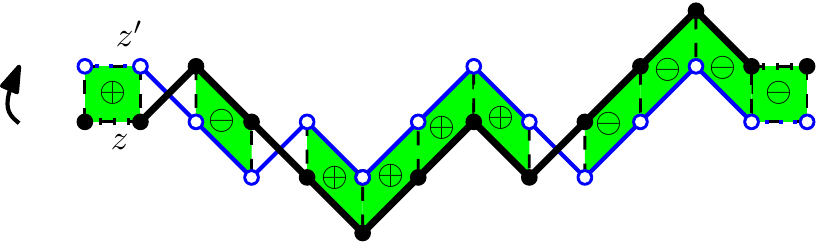}
\end{center}
We need to define $\I^{*}$, a modification of the involution $\I$ from Section~\ref{sec:comb}. Given a marked path $(z;z',p)$, its image $(z;z'',q)$ by $\I^{*}$ is defined as follows. If $p$ is the initial horizontal step, then find the first crossing: $z''$ is obtained by changing $z''$ at this step, which becomes the new mark $q$ (If there is no crossing, so that $z'=z^\pm$, then change it to $z''=z^\mp$ and let $q$ be the final horizontal step). Define $\I^{*}$ symmetrically when $p$ is the final horizontal step. In all other cases, define $\I^*$ as in $\I$, except in the cases where one needs to use the special cases of $L$ and $R$: in this case the mark $q$ will be one of the horizontal steps. It is easily seen that $\I^{*}$ is bijective and reverses signs, so that the first formula of the proposition is proved.

For the second one, we notice as in Section~\ref{sec:comb}  that the l.h.s. can be written as the sum of $sign(z;z',p)$ over all marked paths satisfying $z'\in \Gamma^+(z)$ and $z''\in \Gamma^-(z)$. This happens when $p$ is the initial horizontal step for any $z,z'$ with $z'\in \Gamma^+(z)$, and these cases contribute $3^K$ to the sum. The other possibility is that $p$ is a down step which has no crossing to its left. Clearly $sign(z;z',p)=1$ here also since $z'\in \Gamma^+(z)$, and there are also $3^K$ such possibilities: indeed, their images by $\I^*$ in this case are exactly the marked paths $(z;z'',q)$ where  $z''\in\Gamma^-(z)$ and $q$ is the initial horizontal step. This completes the proof.
\end{proof} 

\begin{figure}[!ht]
\includegraphics[width=0.2\textwidth]{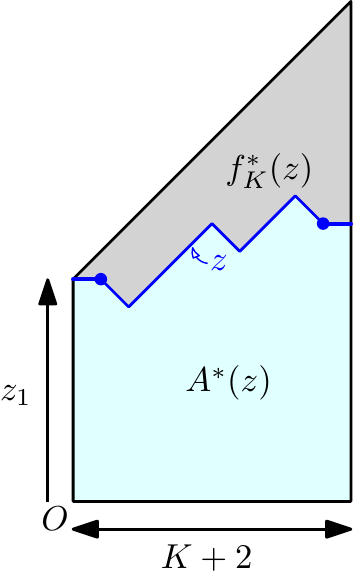}
\caption{Relation between the first coordinate $z_1$ and the area under the path.
\label{Figure_Coupling_2}}
\end{figure}

Denote by $f_K^{*}(z)$ the area between a path $z=(z_1,z_2,\ldots, z_{K},z_{K+1})\in  \mathcal{C}_{K}$,
 and the segment $(z_1, z_1+1, ...., z_1+K+1,z_1+K+2)$.
Then, for any $z\in  \mathcal{C}_{K}$ (see Figure~\ref{Figure_Coupling_2}),
\begin{equation}\label{aire*}
f_{K}^{*}(z) +A^{*}(z) = (K+2) z_1 + \frac{(K+2)^2}{2}.
\end{equation}

Since $f_{K}^{*}$ is bounded by $(K+2)^2$, the position of the first walker and the area under the extended path can be asymptotically related as 
\begin{equation}\label{aire2}
\frac{Z_{[nt],1}^{*}}{\sqrt{n}} = \frac{1}{K+2}\frac{ A^{*} \big(Z_{[nt]}^{*}\big)}{\sqrt{n}}+o(1).
\end{equation}
The first assertion in Proposition~\ref{prop1} implies that the sequence $(A^{*} (Z_n^{*}))_{n\geq 0}$ 
is a martingale with respect to the natural filtration of the Markov chain.
The main result of~\cite{BCEN} is then deduced from the functional central limit theorem for the rescaled martingale $(A^{*}(Z_{n}^{*})/\sqrt{n})_{n\geq 1}.$
The proof is similar to the one given in Section \ref{sec:proba} in the constrained case. The asymptotic variance is then rewritten as
\begin{eqnarray*}
\sigma_{K,*}^{2} &=& \frac{1}{K+2}\  \frac{A_{K}^{*}}{B_{K}^{*}}
\end{eqnarray*}
where
$$A_{K}^{*}= 2 \sum_{\stackrel{\small z,z'\in\mathcal{C}_{K}}{ z'\in\Gamma^+(z)}} (A^{*}(z')-A^{*}(z))\ 
\text{ and }\ 
B_{K}^{*}=2 |\mathcal{M}_{K}|=2\cdot 3^K.$$
The second assertion in Proposition~\ref{prop1} leads to $\sigma_{K,*}^{2} = \frac{2}{K+2}.$ 

\section{Discussion and open problems}
\label{sec:conclusion}

\begin{enumerate}
\item The case when the distance $h$ between the extremal coordinates depends on $K$ can also be considered. It is not difficult to see that Theorem \ref{Th1}
still holds with $h$ replaced by $h(K)$ in the variance formula. The asymptotic behavior of the variance for $K$ large (item i) in Theorem \ref{Th1}) is still true if $h(K)=o(K^{3/4})$. 
Indeed, relation (\ref{tll}) is still valid and local limit theorem 2.3.11 in \cite{LL} (p.46) gives the result after elementary computations. 
The asymptotic behavior of the variance for $h \gg K^{3/4}$ is not known. 

\item Instead of $ \mathcal{C}_{K,h}$ we can consider the following set of paths  
 $$ \mathcal{C}_{K}=\{ z\in \mathbb{Z}^{K+1}; \forall i\in [\![1;K]\!], |z_{i+1}-z_i| \in\{0,1\} \}$$ 
and $(Z_n)_{n\geq 0}$  the Markov chain defined on  $\mathcal{C}_{K}$ corresponding to $K+1$
 simple random walks on $\mathbb{Z}$ under the shape constraint. 
 In \cite{BCEN}, the  set of paths
 $$ \mathcal{D}_{K}=\{ z\in \mathbb{Z}^{K+1}; \forall i\in [\![1;K]\!], |z_{i+1}-z_i| =1 \}$$  
and the corresponding  Markov chain  $(Z_n^{(K)})_{n\geq 0}$ were studied. Whatever the distribution of $Z_0^{(K)}$, the sequence of 
random variables $(Z_{n,1}^{(K)}/\sqrt{n})_n$ converges to a centered gaussian law with variance $2/(K+2)$. By remarking that given the number of zeroes $N_0$ of  the initial random variable $Z_0$,
the position of the first walkers $(Z_{n,1})_{n}$ and $(Z_{n,1}^{(K-N_0)})$ are identically distributed (with the  following convention: if $N_0=K$, the first random walker evolves as the simple symmetric random walk 
on $\mathbb{Z}$), it follows that the sequence of random variables $(Z_{n,1}/\sqrt{n})_n$ converges in law to a mixture of gaussian distributions with variance equal to 
$$\sigma_K^2 = 2 \sum_{l=0}^{K} \frac{  \PP[ N_0 =l]}{K+2-l}.$$
Note that whatever the distribution of the random variable $N_0$ the variance is still greater than the one obtained in the unconstrained case considered in \cite{BCEN} and is rational when $Z_0$ is uniformly distributed on $ \mathcal{C}_{K}$.

\item The generalization of Theorem \ref{Th1} to more general sets of paths does not seem to be obvious.
Except in the case of \cite{BCEN} where our proof can be adapted the construction of a convenient martingale from which the computation of the variance can be deduced is far from clear. 
For instance consider the set of Bernoulli bridges with length $K$ with the additional constraint: $L$ of the $K+1$ random walkers stay at the same height at each step. 
A central limit theorem for the first random walker should also hold. We conjecture that the variance should be increasing in $L$.

\end{enumerate}
       \medskip
    
{\noindent\bf Acknowledgments.} We are grateful to James Norris and Serge Cohen for stimulating discussions.

\end{document}